\documentclass[letterpaper, 10 pt, conference]{ieeeconf}

\IEEEoverridecommandlockouts
\usepackage{color}
\usepackage{amsmath, amssymb}

\usepackage{multicol}
\usepackage{nicefrac}
\usepackage{graphicx}
\usepackage{epstopdf}
\usepackage{tikz}
\usepackage{bm}
\usepackage{picture}
\usepackage{epstopdf}
\usepackage{amsmath}
\usepackage{amsfonts}
\usepackage{amssymb}
\usepackage{cite}

\newcommand{\R}{\mathbb{R}}

\def\QED{\mbox{\rule[0pt]{1.3ex}{1.3ex}}}

\allowdisplaybreaks

\title{Model Predictive Path Planning in Navier–Stokes Flow with POD-Based Reduced-Order Models}
\author{Adam Waterman$^\dagger$, Martin Guay$^\dagger$
\thanks{$^\dagger$Adam Waterman and Martin Guay are with the Department of Chemical Engineering,
		Queen's University, Kingston, ON, K7L 3N6, Canada. {\tt\small 13adw@queensu.ca, guaym@queensu.ca}.}}

\begin{document}
\maketitle

\begin{abstract}
We present a framework for optimal trajectory generation in flow-driven systems governed by the Navier–Stokes equations, combining a Proper Orthogonal Decomposition (POD) reduced-order model (ROM) with Model Predictive Control (MPC). The approach (i) approximates the velocity field from data via snapshot POD and orthogonal projection, (ii) derives a Galerkin-projected dynamical model in reduced coordinates, and (iii) employs MPC to plan control inputs that steer an agent through the predicted flow while satisfying state and actuation constraints. By leveraging reduced-order modeling, the method enables real-time control in high-dimensional flow environments. Simulations demonstrate accurate flow-field reconstruction and efficient trajectory generation within realistic wind environments.

\end{abstract}

\section{Introduction} \label{sec0}

Model Predictive Control (MPC) has become a leading framework for constrained multivariable control in process, power, automotive, and aerospace systems~\cite{Oldewurtel2012,Sturzenegger2016,Noack2003,Bergmann2005}. 
At each sampling instant, MPC solves a finite-horizon optimization problem that balances tracking performance with state and actuator constraints, applying only the first control action before re-optimizing at the next step. 
Despite its wide success, MPC remains computationally demanding when applied to systems described by high-fidelity or infinite-dimensional models, such as those arising from the Navier--Stokes (NS) equations in fluid dynamics. 
Repeatedly solving optimization problems with such models can be prohibitive for real-time implementation, particularly when the prediction horizon must be long enough to capture advective transport phenomena.

Model reduction techniques address this issue by projecting the system dynamics onto a low-dimensional subspace while retaining dominant input--output behavior. Classical linear approaches include balanced truncation and moment-matching methods, while data-driven techniques such as Dynamic Mode Decomposition (DMD) and Proper Orthogonal Decomposition (POD) have proven effective for nonlinear and large-scale systems~\cite{Benner2015,Hinze2009,Antoulas2005,yano2020model}. Comprehensive reviews such as Yano~\cite{yano2020model} highlight the integration of model reduction with optimal control and MPC in aerodynamic and PDE-constrained systems. In fluid mechanics, POD combined with Galerkin projection yields reduced-order models (ROMs) that efficiently capture coherent flow structures and their temporal evolution~\cite{Rowley2005,Noack2011}.

POD was first introduced by Lumley in the context of turbulent flow analysis~\cite{lumley1967structure} and later formalized for reduced-order modeling by Berkooz, Holmes, and Lumley~\cite{berkooz1993proper}. It has since become a cornerstone of reduced-order modeling in fluid mechanics~\cite{Lumley2007,Taira2023}. Using the snapshot method, POD identifies orthonormal spatial modes that maximize the kinetic energy captured from a set of velocity-field snapshots in a least-squares sense.

Truncating the expansion to the most energetic modes produces a low-dimensional basis that accurately reconstructs flow dynamics within the subspace. 
Substituting this representation into the NS equations and applying Galerkin projection yields a nonlinear system of ordinary differential equations (ODEs) for the modal coefficients, which can be estimated using nonlinear filtering methods such as the Extended Kalman Filter (EKF). 
There has been increased interest in the integration of reduced-order models within MPC frameworks (ROM--MPC), spanning building climate control~\cite{Oldewurtel2012,Sturzenegger2016}, flow and aerodynamic load control~\cite{Noack2003,Bergmann2005}, and other PDE-constrained problems~\cite{Benner2015,Hinze2009,Antoulas2005}. 
Many works derive Galerkin ROMs from the incompressible NS equations for controller design~\cite{Rowley2005,Noack2011}, while others identify low-order surrogates directly from data and embed them in MPC formulations~\cite{Korda2018,Proctor2016}. 

Recent studies have also investigated autonomous navigation of
stratospheric balloons using data-driven and learning-based approaches.
The Balloon Learning Environment (BLE)~\cite{greaves2021autonomous,greaves2021ble}
provides a high-fidelity, open-source simulation platform that models realistic atmospheric winds and supports reinforcement-learning-based navigation. In contrast, this work pursues a physics-informed, model-based framework for predictive path planning.

This work extends the POD-based estimation frameworks developed in~\cite{Guay2008,Park2015} to generate optimal trajectories through flow fields predicted by a reduced-order model. 
Rather than controlling a specific dynamic system, the proposed formulation focuses on \emph{path generation}---using predicted velocity fields to plan trajectories that exploit vertical shear or layerwise advection in geophysical flows such as atmospheric wind or ocean currents. 
The framework couples a POD--Galerkin ROM of the NS equations with a receding-horizon optimization problem, enabling computationally efficient, flow-aware trajectory planning.

The remainder of this paper is organized as follows. 
Section~\ref{sec:background} reviews the POD-based model reduction and observer formulation. 
Section~\ref{sec:pathplanning} presents the model predictive path-planning framework that uses the ROM-predicted flow field. 
Section~\ref{sec:simulation} describes the simulation setup and results, and Section~\ref{sec:conclusion} concludes the paper with remarks and directions for future research.
\section{Background}
\label{sec:background}

This section summarizes the Proper Orthogonal Decomposition (POD)–based model reduction and observer design approach used in this work. 
A detailed derivation can be found in~\cite{Guay2008,Park2015}, while comprehensive discussions of POD theory are available in~\cite{Lumley2007,Rowley2005,Noack2011,Taira2023}.

\subsection{POD-Based Model Reduction}

We consider incompressible fluid flow governed by the Navier--Stokes (NS) equations:
\begin{equation} \label{eq:NS_incomp}
\begin{aligned}
div(\mathbf{v}) &= 0 \\
\frac{\partial \mathbf{v}}{\partial t} &= -(\mathbf{v}\cdot\nabla)\mathbf{v} + \nu\nabla^2\mathbf{v} - \nabla p
\end{aligned}
\end{equation}
where $\mathbf{v} : \Omega \times \R \rightarrow \R^3$ indicates velocity field on the spatial domain $\Omega$, $p$ is the pressure term, $\nu = 1/Re$, $Re$ is the Reynolds number. It is assumed that the velocity and pressure field are defined on a closed subset of $\R^3$. Equation \eqref{eq:NS_incomp} represents a scaled formulation of the NS equation where the velocities are scaled by a factor of $V$, time by $V/L$, pressure by $\rho V^2$ where $\rho$ is the density and the viscosity by $\rho VL$, where V and L are nominal velocities and length.

The velocity field is approximated by a truncated series of spatial basis functions (POD modes) 
$\{\boldsymbol{\phi}_{i}(\mathbf{x})\}_{i=1}^{n}$:
\begin{equation}
\mathbf{v}(\mathbf{x},t)
=\sum_{i=1}^{n}a_{i}(t)\,\boldsymbol{\phi}_{i}(\mathbf{x}),
\label{eq:PODexpansion}
\end{equation}
where $a_{i}(t)$ are time-dependent modal coefficients.  
The projection onto the modes requires an inner product over a Hilbert space $\mathcal{H}$ defined as

\begin{equation} \label{HilbProd}
    \langle \mathbf{v}_i, \mathbf{v}_j \rangle = \int_\Omega\mathbf{v}_i(x)\cdot\mathbf{v}_j(x)dV
\end{equation}
where $\langle \mathbf{v}_i, \mathbf{v}_j \rangle$ is the standard inner product between two vectors $\mathbf{v}_i$ and $\mathbf{v}_j$ in Euclidean space, $dV$ is a volume element on $\R^3$. 
The POD based model reduction restates the NS equation using the expression \eqref{eq:PODexpansion}. By substituting ~\eqref{eq:PODexpansion} into~\eqref{eq:NS_incomp}:
\begin{equation}
\begin{aligned}
    \frac{\partial \mathbf{v}}{\partial t} &= \sum_{i=1}^n\dot{a_i}(t)\phi(x) \\
    & = -(\sum_{j=1}^na_j(t)\phi_j(x)\cdot\nabla)\sum_{k=1}^na_k(t)\phi_k(x) \\
    &+ \nu\sum_{i=1}^na_i(t)\nabla^2\phi_i(x) - \nabla p
\end{aligned}
\end{equation}

Projecting onto the space of POD modes $\phi_i(x)$ leads to

\begin{equation} \label{POD_proj}
    \begin{aligned}
       & \langle\sum_{k=1}^n\dot{a}_k(t)\phi_k(x), \phi_i(x)\rangle = \\
       &-\langle\bigg(\sum_{j=1}^n a_j(t)\phi_j(x), \phi_i(x)\cdot\nabla\bigg)\sum_{k=1}^n a_k(t)\phi_k(x),\phi_i(x) \rangle \\
       &+ \nu\langle\sum_{i=k}^na_k(t)\nabla^2\phi_k(x), \phi_i(x)\rangle - \langle\nabla p,\phi_i(x)\rangle
    \end{aligned}
\end{equation}

By orthogonality of the modes

\begin{equation}
\langle\phi_i , \phi_j(x) \rangle = 
    \begin{cases}
    1, & i = j, \\
    0, & i\neq j.
\end{cases}
\end{equation}
For divergence-free modes with $\boldsymbol{\phi}_{i}\!=\!0$ on $\partial\Omega$, 
the pressure term vanishes under integration by parts.  
The reduced-order model (ROM) \eqref{POD_proj} therefore becomes

\begin{equation} \label{Coeff_dyn}
\begin{aligned}
   \dot{a}_i(t)&\langle\phi_i(x), \phi_k(x)\rangle = \dot{a}_k(t) = \\ &-\sum_{i=1}^n\sum_{j=1}^na_i(t)a_j(t)\langle(\phi_j(x)\cdot\nabla)\phi_i(x),\phi_k(x) \rangle \\
   &+ \nu\sum_{i = 1}^na_i(t)\langle\nabla^2\phi_i(x),\phi_k(x)\rangle - \langle\nabla p, \phi_k(x)\rangle.
\end{aligned}
\end{equation}

Equation \eqref{Coeff_dyn} is the decomposition model of the velocity field $v(x,t)$ on the POD modes. This finite dimensional approximation provides the basis for the design of an observer for the fluid flow dynamic system. 

\subsection{Computation of POD Modes}

The POD basis functions are derived from a representative ensemble of $m$ velocity snapshots 
$\{\mathbf{u}_{k}(\mathbf{x})\}_{k=1}^{m}$ obtained from simulation or experiment.  
The goal is to find an orthonormal set $\boldsymbol{\phi}_{i}$ that minimizes the mean-square reconstruction error
\begin{equation}
\min_{\boldsymbol{\phi}_{1},\ldots,\boldsymbol{\phi}_{n}}
\frac{1}{m}\!\sum_{k=1}^{m}
\left\|
\mathbf{u}_{k}
-\!\!\sum_{i=1}^{n}
  \langle \mathbf{u}_{k},\boldsymbol{\phi}_{i}\rangle
  \boldsymbol{\phi}_{i}
\right\|^{2}.
\label{eq:PODopt}
\end{equation}
Solving~\eqref{eq:PODopt} leads to the eigenvalue problem
\begin{equation}
\mathbf{U}\mathbf{c}=\lambda\mathbf{c},
\quad
U_{ij}=\tfrac{1}{m}\langle\mathbf{u}_{i},\mathbf{u}_{j}\rangle,
\label{eq:eig}
\end{equation}
where $\lambda$ represents the energy associated with each mode.  
The resulting modes are linear combinations of snapshots,
\begin{equation}
\boldsymbol{\phi}_{\ell}
=\sum_{k=1}^{m}c_{k\ell}\,\mathbf{u}_{k},
\label{eq:modes}
\end{equation}
ordered by decreasing $\lambda_{\ell}$ so that the first few modes capture the majority of the kinetic-energy content~\cite{Lumley2007,Rowley2005,Taira2023}.  

\subsection{Observer Design}

The dynamical system \eqref{Coeff_dyn} yields a set of quadratic differential equations of the form:
\begin{equation} \label{ROMDyn}
    \dot{a}_k(t) = L_ka(t) + a(t)^T\mathcal{Q}_ka(t)
\end{equation}
where $L_k$ is a row vector with elements given by 

\begin{equation}
    L_{ik} = \langle\nabla^2\phi_i(x), \phi_k(x)\rangle
\end{equation}
and $\mathcal{Q}_k$ is an $n$ by $n$ matrix with elements

\begin{equation}
    \mathcal{Q}_ijk = \langle(\phi_j(x)\cdot\nabla)\phi_i(x), \phi_k(x)\rangle
\end{equation}
for $k = 1,...,n$. 

Since the POD modes are such that $div(\phi) = 0$, it follows that:

\begin{equation}
    \begin{aligned}
        \int_\Omega \nabla p \cdot\phi_k(x)dV &= \int_\Omega div(p\phi_k(x))dV \\
        &= \int_{\partial \Omega}p\phi_k(x)\cdot\mathbf{n}_\Omega dS
    \end{aligned}
\end{equation}
where $\mathbf{n}_\Omega$ represents the unit vector normal to the spatial domain $\Omega$. Therefore the pressure term will vanish over a closed domain $(\phi_k(x) = 0$) on the boundary of $\Omega$, $\partial \Omega$ and can be ignored.

It is assumed that several velocity measurements $v_0$ are available at predefined locations $x_0$. The measurements can be expressed using the expansion \eqref{eq:PODexpansion}. For example, if one measures the average velocity $\mathbf{v}_{avg}(x_0,t) = (\mathbf{v}_1(x_0,t) + \mathbf{v}_2(x_0,t) + \mathbf{v}_3(x_0,t))$, at a point $x_0$, then the corresponding measurement becomes

\begin{equation}
    \mathbf{v}_{avg}(x_0, t) = \sum_{i=1}^na_i(t)(\phi_i^{(1)}(x) + \phi_i^{(2)}(x) + \phi_i^{(3)}(x))  
\end{equation}
where $\phi_i^j$ represents the $j^{th}$ element of the $i^{th}$ POD mode. Since the POD modes are time independent, the resulting output relation can
be written as
\begin{equation}
    y_{\text{fixed}}(t) = C_{\text{fixed}} a(t),
\end{equation}
where $C_{\text{fixed}}\in\mathbb{R}^{m\times n}$ is the measurement matrix
obtained by evaluating the POD modes at $m$ fixed sensor locations, and
$a(t)\in\mathbb{R}^n$ is the vector of time-varying Galerkin coefficients
of the reduced-order velocity field $\mathbf{v}(\mathbf{x},t)$.

In addition to the fixed sensor network, the navigating agent is assumed to
act as a mobile sensor. The corresponding time-varying measurement output is
\begin{equation}
    y_{\text{mobile}}(t) = C_{\text{mobile}}(t)a(t),
\end{equation}
where $C_{\text{mobile}}(t)$ represents the POD modes evaluated at the
agent's instantaneous position. The complete measurement vector is formed by
stacking the fixed and mobile measurements:
\begin{equation}
    y(t) =
    \begin{bmatrix}
        y_{\text{fixed}}(t) \\[2mm]
        y_{\text{mobile}}(t)
    \end{bmatrix},
    \qquad
    C(t) =
    \begin{bmatrix}
        C_{\text{fixed}} \\[2mm]
        C_{\text{mobile}}(t)
    \end{bmatrix},
\end{equation}
so that the overall measurement equation becomes
\begin{equation}
    y(t) = C(t)a(t).
\end{equation}

The complete reduced-order system with available measurements is therefore
\begin{equation} \label{eq:meas}
    \begin{aligned}
        \dot{a}_k(t) &= L_ka(t) + a(t)^T\mathcal{Q}_ka(t), \quad k = 1,...n \\
        y(t) &= C(t)a(t).
    \end{aligned}
\end{equation}

where $L$ and $\mathcal{Q}$ are defined by the Galerkin projection in
\eqref{ROMDyn}. The use of a time-varying measurement matrix $C(t)$
reflects the changing sensing geometry as the agent moves through the flow
field. Each row of $C(t)$ corresponds to evaluating the POD modes
$\{\boldsymbol{\phi}_i(\mathbf{x})\}$ at the current sensor positions
$\mathbf{x}_s(t)$. When all sensors are stationary, $C(t)$ reduces to a
constant matrix $C_{\text{fixed}}$. This formulation allows the observer to
assimilate measurements from mobile platforms and continuously update the
state estimate $\hat{\mathbf{a}}(t)$ based on a time-dependent measurement
map, thereby improving flow-field reconstruction in sparsely instrumented
environments.

An Extended Kalman Filter (EKF) is used to estimate the reduced--order
state $\hat{\mathbf{a}}(t)$ from the nonlinear system~\eqref{eq:meas}.
At each step, the nonlinear ROM dynamics and the time--varying measurement
map are linearized about the current estimate to obtain
\begin{align}
    \dot{\hat{\mathbf{a}}}(t)
    &= L\hat{\mathbf{a}}(t)
       + \mathcal{Q}(\hat{\mathbf{a}}(t))
       + K(t)\!\left[y(t) - C(t)\hat{\mathbf{a}}(t)\right], \label{eq:EKFstate}\\[1mm]
    \dot{P}(t)
    &= A(t)P(t) + P(t)A(t)^{\!\top}
       - P(t)C(t)^{\!\top}R^{-1}C(t)P(t) + Q, \label{eq:EKFcov}
\end{align}
where $P(t)$ is the state--error covariance, $Q$ and $R$ are the process and
measurement noise covariances, and
\begin{align}
    A(t) &= \left.\frac{\partial}{\partial \mathbf{a}}
            \!\left[L\mathbf{a} + \mathcal{Q}(\mathbf{a})\right]\!
            \right|_{\mathbf{a}=\hat{\mathbf{a}}(t)}, \\[1mm]
    K(t) &= P(t)C(t)^{\!\top}\!\left[C(t)P(t)C(t)^{\!\top}\!+\!R\right]^{-1}
\end{align}
is the Kalman gain.  The term $\mathcal{Q}(\mathbf{a})$ represents the
quadratic nonlinear interactions of the ROM dynamics, whose Jacobian
$A(t)$ is obtained analytically from~\eqref{ROMDyn}.  The EKF thus propagates the reduced--order
flow model forward in time while correcting the estimated coefficients
using both fixed and mobile sensor measurements through the
time--dependent matrix $C(t)$.

For numerical implementation, the EKF is applied in discrete time with
sampling period $\Delta t$. The reduced--order model is integrated over
each interval to yield the discrete prediction map
\[
\mathbf{a}_{k+1} = f_d(\mathbf{a}_k) + w_k,
\qquad
y_k = C_k \mathbf{a}_k + v_k,
\]
where $f_d(\cdot)$ represents the ROM dynamics integrated over $\Delta t$,
and $w_k$ and $v_k$ denote process and measurement noise.
The discrete--time EKF proceeds as follows:

\paragraph*{Prediction step}
\begin{align}
\hat{\mathbf{a}}_{k|k-1} &= f_d(\hat{\mathbf{a}}_{k-1}), \\
P_{k|k-1} &= A_{k-1} P_{k-1|k-1} A_{k-1}^{\!\top} + Q,
\end{align}

\paragraph*{Update step}
\begin{align}
K_k &= P_{k|k-1} C_k^{\!\top}
       \left(C_k P_{k|k-1} C_k^{\!\top} + R\right)^{-1}, \\
\hat{\mathbf{a}}_{k|k} &= \hat{\mathbf{a}}_{k|k-1}
       + K_k\!\left(y_k - C_k\hat{\mathbf{a}}_{k|k-1}\right), \\
P_{k|k} &= (I - K_k C_k) P_{k|k-1}.
\end{align}

Here $A_{k-1}$ is the discrete--time Jacobian of $f_d(\cdot)$ evaluated
at $\hat{\mathbf{a}}_{k-1}$, which can be obtained analytically from the
ROM coefficients or approximated numerically via finite differences.
This formulation provides a practical recursive estimator compatible with
the sampling rate used in the ROM--based wind prediction and
path--planning simulations.

The estimated coefficients can then reconstruct the full velocity field via~\eqref{eq:PODexpansion}, providing an updated prediction 
\begin{equation}
v_{\mathrm{pred}}(\mathbf{x},t)
=\sum_{i=1}^{n}\hat{a}_{i}(t)\,\boldsymbol{\phi}_{i}(\mathbf{x}),
\label{eq:vpred}
\end{equation}
which serves as the forecasted flow field for the predictive path-planning framework introduced in Section~\ref{sec:pathplanning}.

%

\section{Model Predictive Path Planning Using a ROM-Predicted Flow Field}
\label{sec:pathplanning}

\subsection{Predicted Flow Environment}

The reduced-order model (ROM) developed in Section~\ref{sec:background} provides a compact representation of the flow field governed by the Navier--Stokes equations. 
The ROM state $\mathbf{a}(t)\in\mathbb{R}^{r}$ evolves autonomously from \eqref{ROMDyn}

\begin{equation*}
    \dot{a}_k(t) = L_ka(t) + a(t)^TQ_ka(t)
\end{equation*}
and reconstructs the velocity field through \eqref{eq:vpred}.

Equation~\eqref{eq:vpred} defines a time-varying velocity field predicted from the ROM coefficients, representing the ambient flow that will advect a passive or lightly actuated agent.

In many flow-driven systems---such as high-altitude balloons or underwater gliders---horizontal motion is dominated by advection from the surrounding flow, while actuation is primarily available in the vertical direction. 
For demonstration, the agent is treated as a point mass advected by the predicted flow whose kinematics follow
\begin{equation}
\dot{\mathbf{x}}(t)
=\mathbf{v}_{\mathrm{pred}}\big(\mathbf{x}(t),t\big)
+u_{z}(t)\,\mathbf{e}_{z},
\label{eq:agentdyn}
\end{equation}
where $\mathbf{x}=[x,\,y,\,z]^{\top}\!\in\!\mathcal{X}\!\subset\!\mathbb{R}^{3}$ is the agent position, 
$u_{z}(t)\!\in\!\mathcal{U}\!\subset\!\mathbb{R}$ is the controllable vertical velocity, and $\mathbf{e}_{z}=[0,0,1]^{\top}$ denotes the unit vector in the vertical direction.

\subsection{Path Generation via Model Predictive Optimization}

At discrete time $t_{k}$, the ROM coefficients $\hat{\mathbf{a}}_{k}$ provide a forecast of the velocity field $\mathbf{v}_{\mathrm{pred}}(\mathbf{x},t)$ over a finite horizon $T=N\Delta t$. 
The goal of the planner is to determine an optimal sequence of vertical inputs 
$\{u_{z,k+i|k}\}_{i=0}^{N-1}$ that steers the agent toward a desired target $\mathbf{x}_{\mathrm{ref}}$.

The discrete-time prediction model derived from~\eqref{eq:agentdyn} is
\begin{equation}
\mathbf{x}_{k+i+1|k}
=\mathbf{x}_{k+i|k}
+\Delta t\,
\Big[
v_{\mathrm{pred}}\!\big(\mathbf{x}_{k+i|k},t_{k+i}\big)
+u_{z,k+i|k}\mathbf{e}_{z}
\Big],
\label{eq:discrete}
\end{equation}
for $i=0,\ldots,N-1$.
The corresponding finite-horizon optimization problem is formulated as
\begin{align}
\min_{\{u_{z,k+i|k}\}} & \quad
V_{f}\big(\mathbf{x}_{k+N|k}\big)
+\sum_{i=0}^{N-1}\!
\ell\big(\mathbf{x}_{k+i|k},u_{z,k+i|k}\big)
\nonumber\\[2mm]
\text{s.t.} \label{eq:MPCopt} \\
 \mathbf{x}_{k+i+1|k}
 &=
\mathbf{x}_{k+i|k}+\Delta t\,
\big[
\mathbf{v}_{\mathrm{pred}}(\mathbf{x}_{k+i|k},t_{k+i})
+u_{z,k+i|k}\mathbf{e}_{z}
\big],
\nonumber\\
 \mathbf{x}_{k+i|k}&\in\mathcal{X},\;
u_{z,k+i|k}\in\mathcal{U},\;
\mathbf{x}_{k+N|k}\in\mathcal{X}_{f}.
\nonumber
\end{align}

The stage cost $\ell(\mathbf{x},u_{z})$ penalizes both target deviation and control effort:
\begin{equation}
\ell(\mathbf{x},u_{z})
=w_{p}\,\|\mathbf{x}-\mathbf{x}_{\mathrm{ref}}\|_{2}^{2}
+w_{u}\,u_{z}^{2},
\label{eq:stagecost}
\end{equation}
and the terminal cost $V_{f}$ encourages convergence to the target region.
The feasible sets $\mathcal{X}$ and $\mathcal{U}$ impose spatial and actuation limits, such as altitude bounds or maximum climb rate.

Once~\eqref{eq:MPCopt} is solved, only the first input $u_{z,k|k}^{\star}$ is applied, and the agent state is advanced according to~\eqref{eq:discrete}. 
At the next step, new observations update the ROM coefficients via the observer in Section~\ref{sec:background}, producing an updated velocity forecast for replanning. 
This receding-horizon process yields an adaptive \emph{path generator} that continuously refines the trajectory as the predicted flow evolves.

\subsection{Planning Horizon and Feasibility}

The prediction horizon $T=N\Delta t$ and sampling time $\Delta t$ determine both the planner’s foresight and computational cost. 
A short horizon may cause myopic behavior, whereas an excessively long horizon can make the nonlinear optimization intractable or degrade performance by overly trusting longer open-loop predictions 
In practice, these parameters are selected based on the characteristic advection timescale of the flow and computational capability. 
Real-time feasibility is achieved using Sequential Quadratic Programming or real-time iteration (RTI) schemes~\cite{Diehl2005}.

Because the agent’s horizontal motion is entirely driven by the predicted NS-based flow field, strict recursive feasibility or cost-decrease guarantees cannot generally be ensured. 
Instead, the planner acts as a best-effort trajectory generator: at each step it computes a feasible, energy-efficient path under the current ROM forecast and re-optimizes as new information becomes available. 
This formulation allows the same predictive strategy to be applied to a wide range of flow-driven platforms without modification to the underlying dynamics. Unlike traditional ROM--MPC formulations designed for closed-loop stabilization~\cite{Benner2015,Hinze2009,Antoulas2005}, 
the present work emphasizes predictive path generation and navigation performance in dynamically evolving flow environments.

\section{Simulation Framework and Results}
\label{sec:simulation}

\subsection{Overview}

To demonstrate the proposed model‐predictive path‐planning framework, a simulation study was performed using a reduced‐order model (ROM) of a three‐dimensional flow field derived from the incompressible Navier–Stokes equations.  
The objective is to evaluate the planner’s ability to generate feasible and efficient trajectories through a time‐varying velocity field predicted by the ROM.  
The simulation illustrates how predicted flow information can be exploited to navigate a passive or lightly actuated agent---represented here as a point mass with altitude‐only control---toward a desired target region. It is assumed that the flow field is adequately observed via a distributed sensor network ensuring observability of the training dataset. There were 8 sensors including the agent used for this simulation. 

The simulation setup parallels the station-keeping problem formulated
in the Balloon Learning Environment (BLE)~\cite{greaves2021autonomous},
where the objective is to maintain a high-altitude balloon within a
specified radius of a target location despite spatiotemporal variations
in wind velocity.
Following this framework, a circular target zone of $50~\mathrm{km}$
radius is defined around the desired station point. The ROM-predicted
velocity field provides the flow advection dynamics, and the
model-predictive planner computes vertical-speed commands to maintain
proximity to the station while minimizing control effort.

\subsection{Reduced‐Order Flow Model}

The flow environment was obtained from the European Centre for Medium-Range Weather Forecasts Re-analysis (ERA5) dataset over a limited spatial domain $\Omega=[x_{\min},x_{\max}]\!\times\![y_{\min},y_{\max}]\!\times\![z_{\min},z_{\max}]$. For the simulation, the domain chosen had a 2 degree range in the latitude and longitude (approx. 146 km x 220 km  square area) centered on $49^{\circ}$N and $-81^{\circ}$W and from pressure levels 125 hPA to 10 hPA (approx. from 15 km to 30 km altitude).
Velocity snapshots were collected from randomly selected days at uniform 1 hour temporal intervals over 6 months representing 36 days of data and processed using the snapshot POD procedure described in Section~\ref{sec:background}. Snapshots were pre-processed to remove mean bias so that the ROM captured the dynamics of the wind field.

\begin{figure}[h]
  \centering
  
  \includegraphics[width=0.9\columnwidth]{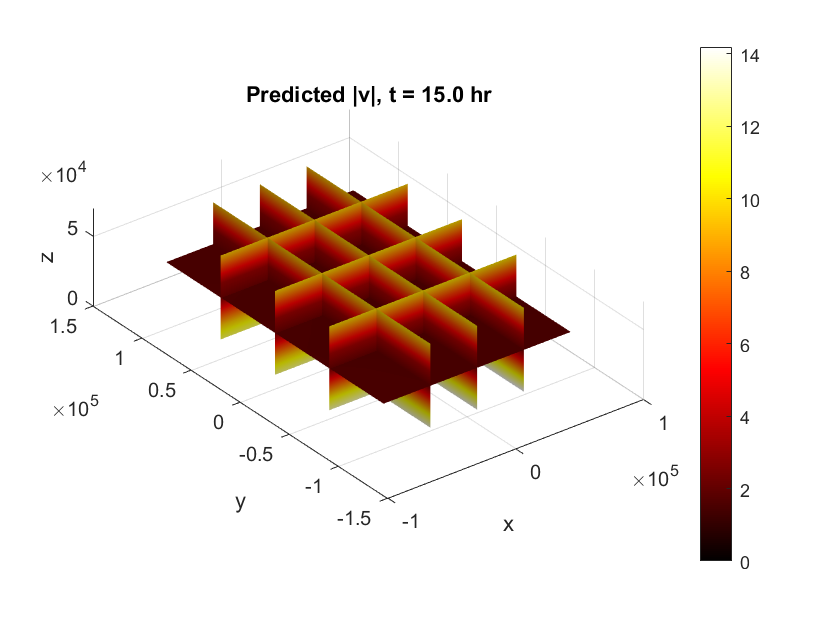}
  \caption{Simulation results showing slices of the predicted mean-squared velocity field.}
  \label{fig:heat_comp_pred}
  
  \vspace{1em} 
  
  \includegraphics[width=0.9\columnwidth]{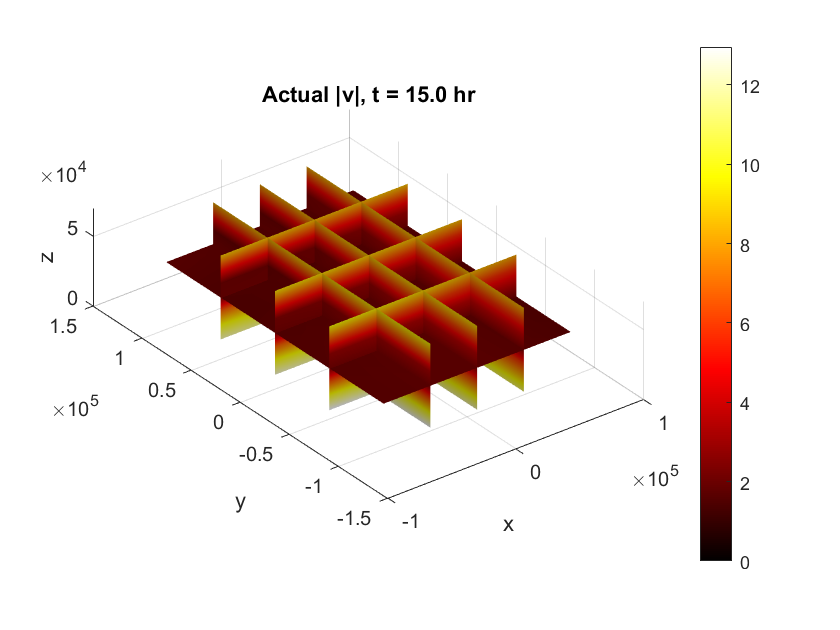}
  \caption{Slices of the actual mean-squared velocity field.}
  \label{fig:heat_comp_actual}
  
\end{figure}

\begin{figure}[h]
  \centering
  \includegraphics[width=0.95\columnwidth]{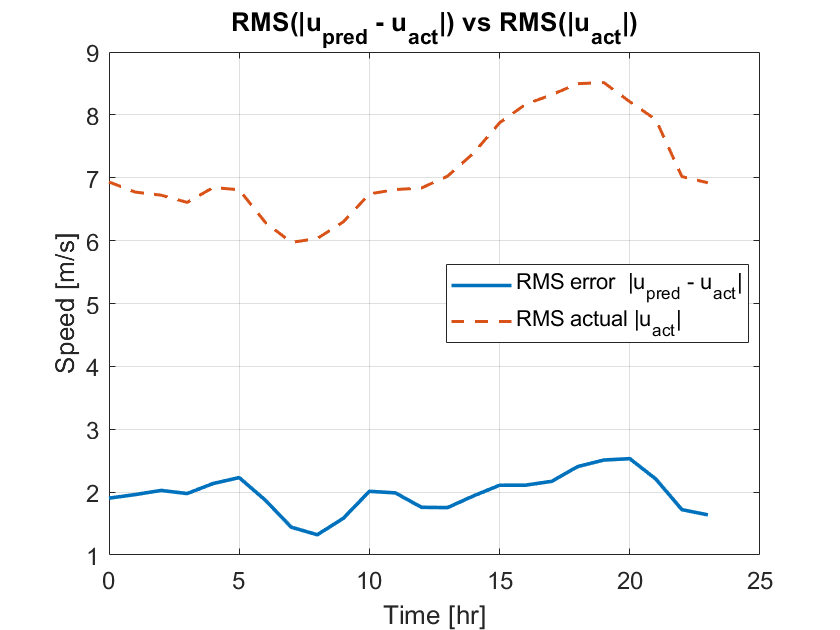}
  \caption{Total root mean squared velocity error between the estimated and actual velocity fields.}
  \label{fig:RMS}
\end{figure}

Figure~\ref{fig:heat_comp_pred} and Figure~\ref{fig:heat_comp_actual} shows the heat map of the predicted and actual mean squared velocity fields for the randomly selected test day at hour 15. Figure~\ref{fig:RMS} shows the average RMSE in the velocity in the estimate of the wind field at each timestep, illustrating the ability of the ROM to predict the actual wind field.

\subsection{Path‐Planning Configuration}

The computational domain was discretized with uniform spatial resolution $\Delta x=\Delta y=\Delta z$ sufficient to resolve dominant flow structures.  
The agent was initialized at $\mathbf{x}_{0}=[x_{0},y_{0},z_{0}]^{\top}$ and tasked to reach a specified target $[x_{ref},y_{ref}]$ within the altitude bounds.
The optimization problem~\eqref{eq:MPCopt} was implemented with horizon length $N$ and sampling time $\Delta t$, such that the prediction horizon $T=N\Delta t$ covered several characteristic advection times of the flow.

The cost‐function weights $w_{p}$ and $w_{u}$ in~\eqref{eq:stagecost} were tuned empirically to balance target proximity and actuation effort.  
Boundary limits represented by $\Omega$ and maximum climb rate $|u_{z}|\!\le\!u_{\max} = 1 $m/s defined the feasible set $\mathcal{U}$.  
At each timestep, the nonlinear optimization was solved via Sequential Quadratic Programming (SQP) with warm starting the initialization from the previous sequence of inputs.

\subsection{Performance Metrics}

The generated trajectories were evaluated using the following metrics:

\begin{enumerate}
\item \textbf{Final target distance}
\begin{equation*}
d_{f}=\|\mathbf{x}(T_{f})-\mathbf{x}_{\mathrm{ref}}\|_{2},
\label{eq:df}
\end{equation*}
where $T_{f}$ is the total simulation time.

\item \textbf{Average alignment with predicted flow}
\begin{equation*}
\bar{\gamma}
=\frac{1}{T_{f}}\!
 \int_{0}^{T_{f}}\!
 \frac{\dot{\mathbf{x}}(t)\!\cdot\!\mathbf{v}_{\mathrm{pred}}(\mathbf{x},t)}
      {\|\dot{\mathbf{x}}(t)\|\,\|\mathbf{v}_{\mathrm{pred}}(\mathbf{x},t)\|}\,dt,
\label{eq:gamma}
\end{equation*}
measuring how effectively the planner exploits favorable flow directions.

\item \textbf{Control effort}
\begin{equation*}
u_{\mathrm{rms}}
=\sqrt{\frac{1}{T_{f}}\int_{0}^{T_{f}}u_{z}^{2}(t)\,dt},
\label{eq:urms}
\end{equation*}
representing average actuation magnitude. 

\item \textbf{Runtime (RT)}

Represented by the total time to run the simulation to demonstrate the computational effort, calculated with MATLAB's tic-toc function. 
\end{enumerate}
Simulations were performed in MATLAB R2024a on a 6-core laptop (AMD Ryzen 2600X, 3.6 GHz, 32 GB RAM).
\subsection{Representative Results}

\begin{figure}[h]
  \centering
  
  \includegraphics[width=1\columnwidth]{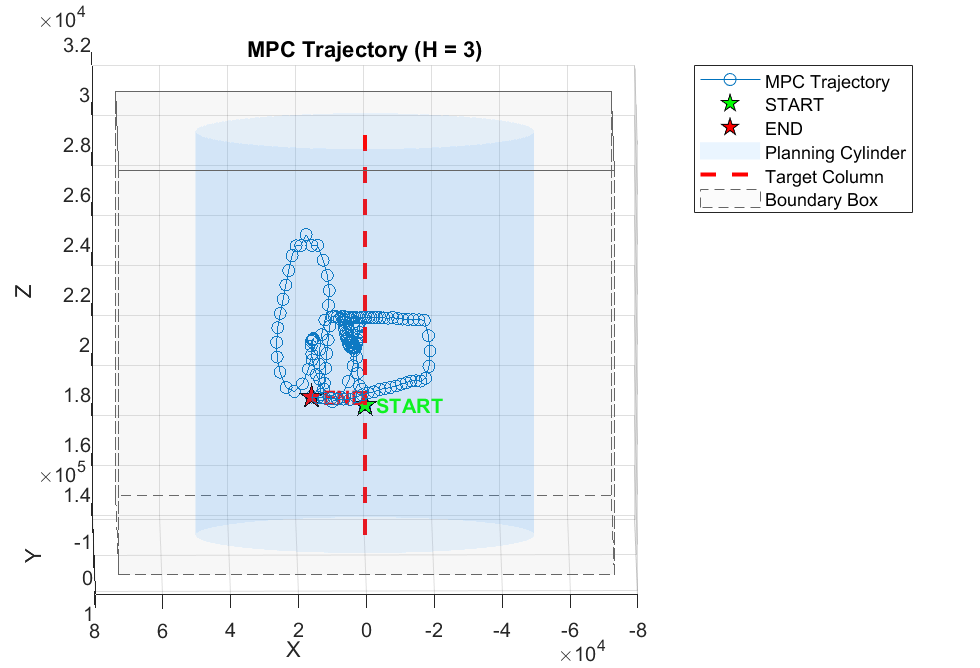}
  \caption{Path generated by the predictive planner with a 3-hour prediction horizon. 
  The grey box indicates the spatial boundary, the blue cylinder is the 50 km radius, and the red dashed marker denotes the target location.}
  \label{fig:traj3hr}
  
  \vspace{1em} 
  
  \includegraphics[width=1\columnwidth]{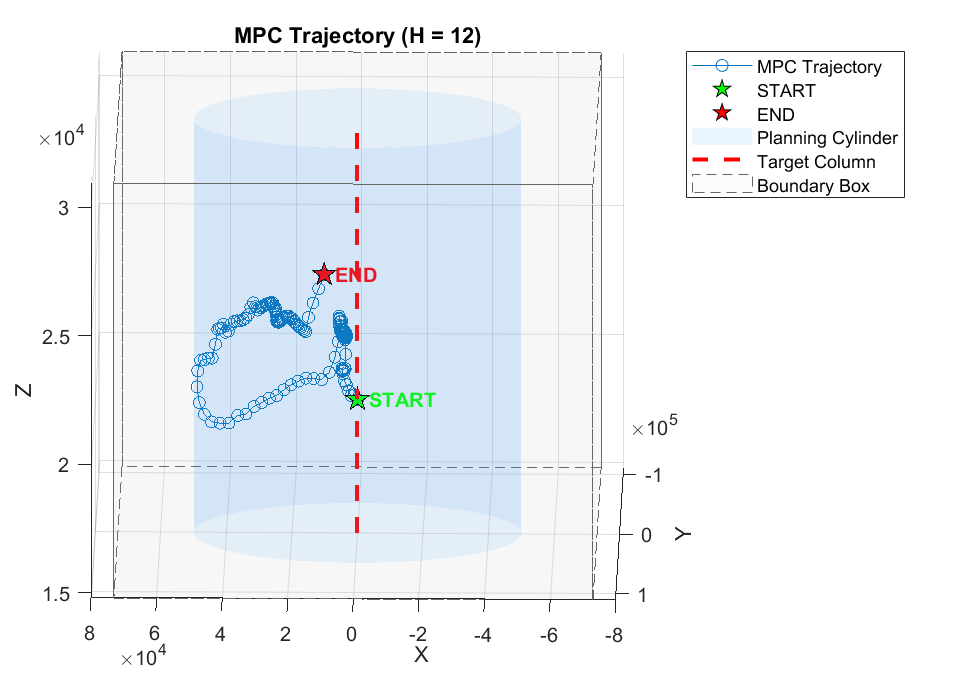}
  \caption{Path generated by the predictive planner with a 12-hour prediction horizon.}
  \label{fig:traj12hr}
  
\end{figure}

\begin{figure}[h]
  \centering
  
  \includegraphics[width=1\columnwidth]{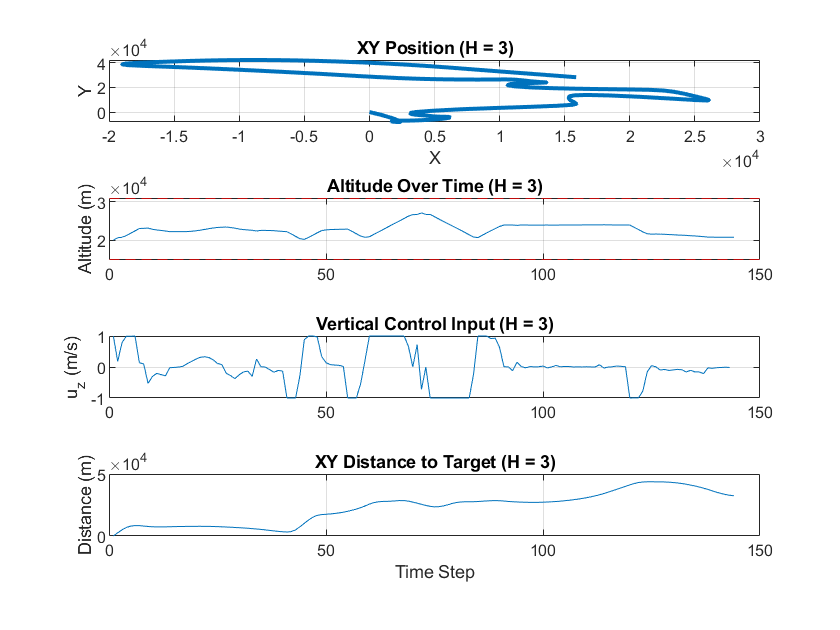}
  \caption{Performance characteristics corresponding to Fig.~\ref{fig:traj3hr}. 
  Altitude adjustments allow the agent to exploit favorable wind strata.}
  \label{fig:3hr_perf}
  
  \vspace{1em} 
  
  \includegraphics[width=1\columnwidth]{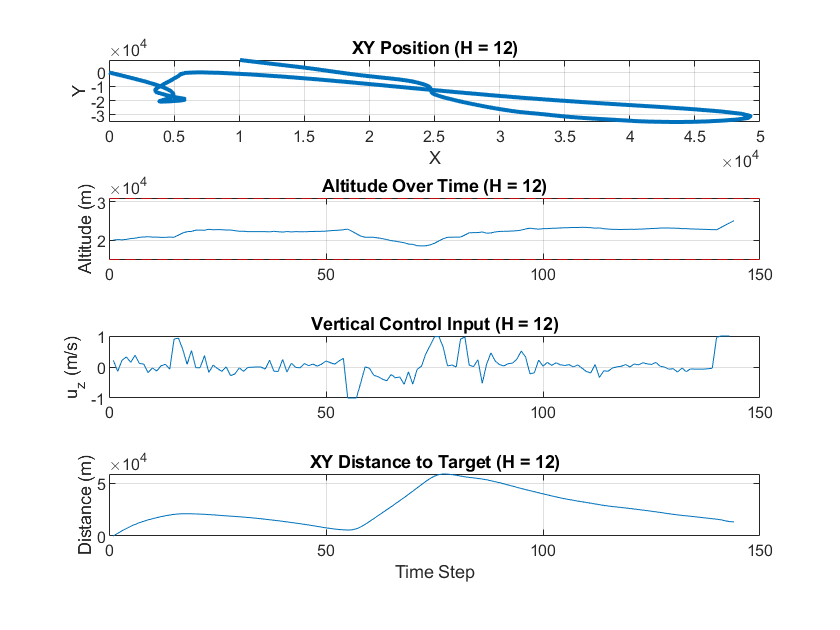}
  \caption{Performance characteristics corresponding to Fig.~\ref{fig:traj12hr}.}
  \label{fig:12hr_perf}
  
\end{figure}

Figure~\ref{fig:traj3hr} and Figure~\ref{fig:traj12hr} shows a representative trajectory produced by the proposed planner for a $3$-hour horizon ($T\!=\!3\,\mathrm{h}$) and $12$-hour horizon ($T\!=\!12\,\mathrm{h}$) respectively with ten-minute sampling intervals over a total 24-hour period.  
The agent begins in a region of west–east advection and adjusts altitude to exploit favorable layers, leading to steady convergence toward the target. The longer prediction horizon can be seen to be more efficient in control action and effective in driving towards the target location as it can avoid short-term gains that cause long-term problems. 
The corresponding performance plots for the 3-hour prediction horizon (Fig.~\ref{fig:3hr_perf}) and the 12-hour prediction horizon (Fig.~\ref{fig:12hr_perf}) illustrates how the optimizer commands vertical maneuvers to transition between flow strata. 

Quantitative performance metrics are summarized in Table~\ref{tab:metrics}, with all metrics averaged over a 24-hour simulation using 10-min sampling. Short horizons produce locally feasible but globally inefficient motion, while excessively long horizons yield limited improvement at higher computational cost.

\begin{table}[h]
\caption{Performance metrics for varying prediction horizons.}
\label{tab:metrics}
\centering
\begin{tabular}{ccccc}
\hline
Horizon (h) & $d_{f}$ (km) & $\bar{\gamma}$ & $u_{\mathrm{rms}}$ (m/s) &  RT (s)\\
\hline
1.5 & 54.6 & 0.92 & 0.60 & 67.0 \\
3.0 & 32.4 & 0.96 & 0.56 & 216.7 \\
6.0 & 13.5 & 0.96 & 0.51 & 696.8 \\
12.0 & 13.4 & 0.96 & 0.36 & 1982.9 \\
\hline
\end{tabular}
\end{table}

\subsection{Discussion}

The results show that coupling POD‐based flow prediction with receding‐horizon optimization enables physically consistent and computationally efficient path generation in complex flow environments.  
The planner adapts to spatiotemporal variations in the predicted velocity field and can be reinitialized as new data becomes available, providing resilience to forecast uncertainty.  
Because the formulation is purely kinematic, it can be integrated with diverse vehicle models as a high‐level guidance module.  
These findings are consistent with prior ROM studies that emphasize the efficiency of reduced‐order modeling for flow prediction~\cite{Rowley2005,Noack2011,Taira2023} and with practical MPC solvers for real‐time optimization~\cite{Diehl2005}.

\section{Conclusion}
\label{sec:conclusion}

This work presented a model-predictive path-planning framework that couples a Proper Orthogonal Decomposition (POD) reduced-order model of the Navier--Stokes equations with predictive optimization to generate flow-aware trajectories. Using a data-driven ROM of the dominant flow structures, the planner treats the agent as a lightly actuated particle navigating predicted wind fields and exploits vertical shear for station keeping and target tracking. Simulations with ERA5 wind data demonstrated that the approach yields physically consistent and computationally efficient paths. Future efforts will extend the framework to incorporate detailed vehicle dynamics and investigate integration with learning-based frameworks such as the Balloon Learning Environment~\cite{greaves2021autonomous,greaves2021ble} to generate a physics-informed, reduced-order hybrid alternative.

\bibliographystyle{IEEEtran}
\bibliography{ROM_MPC_refs}
\end{document}